\documentclass{article}

\usepackage{arxiv}

\usepackage[utf8]{inputenc} 
\usepackage[T1]{fontenc}    
\usepackage{hyperref}       
\usepackage{url}            
\usepackage{booktabs}       
\usepackage{amsfonts}       
\usepackage{nicefrac}       
\usepackage{microtype}      
\usepackage{lipsum}
\usepackage{graphicx}
\graphicspath{ {./images/} }
\usepackage{tabu}
\usepackage{amsmath}
\usepackage{multirow}
\usepackage{threeparttable}
\usepackage{algorithm}
\usepackage{algorithmic}
\usepackage[title]{appendix}

\title{A fast-solved model for energy-efficient train control based on convex optimization}

\author{
 Minling Feng \\
  Shien-ming Wu School of Intelligent Engineering\\
  South China University of Technology\\
  Guangzhou, China 511442 \\
  \texttt{Minling.Feng@outlook.com} \\
   \And
 Kunpeng Wu \\
  Shien-ming Wu School of Intelligent Engineering\\
  South China University of Technology\\
  Guangzhou, China 511442 \\
  \texttt{Kunpeng.Wu@outlook.com} \\
  \And
 Shaofeng Lu* \\
  Shien-ming Wu School of Intelligent Engineering\\
  South China University of Technology\\
  Guangzhou, China 511442 \\
  \texttt{lushaofeng@scut.edu.cn} \\
  }

\begin{document}
\maketitle
\begin{abstract}
In modern rail transportation, energy-efficient train control (EETC) is concerned with the optimal train speed trajectory or control strategies to achieve the minimum energy cost under various operation and traction constraints. This paper proposes an EETC model based on convex optimization so that the model can be rapidly solved by convex optimization algorithms. The high computational efficiency and robustness of the convex model can be verified by comparing the results achieved by the method proposed by this paper and other mainstream mathematical programming methods including mixed-integer linear programming (MILP) and Radau pseudospectral method (RPM). Based on the characteristics of convex optimization, the proposed method boasts more significant advantages over its counterparts in terms of computational efficiency in the promising online applications for automatic train control systems of various types of rail transportation. 
\end{abstract}


\section{Introduction}
\label{sec:introduction}
In the past few decades,  energy-efficient train control (EETC) has been intensively studied with many effective methods developed. Benefitting from the development of computation technology, studies on EETC not only theoretically adopts the optimal control theory to achieve the optimal control strategies but also applies various kinds of direct optimization methods which rely on numeric iteration to seek optimal control strategies of the train.

The classical optimal control theory based on Pontryagin's maximum principle (PMP) was applied to indirectly derive the optimal train control strategies consisting of traction with maximum power, cruising, coasting, and braking with maximum power by analyzing the necessary conditions for optimality. Early works conducted by a number of researchers have successfully applied PMP in locating the optimal train control strategies \cite{Khmelnitsky2000,Howlett2000,Liu2003,Albrecht2016a,Albrecht2016b}. Heuristic algorithm is capable of reducing search space and model complexity with effective heuristic rules to address the EETC problem. Such methods usually take a relatively longer time to search for the optimal solution and thus suffer from the computational burden. Due to the fact that EETC itself alone can be solved by other efficient methods, heuristic methods are more commonly applied to the complex problem of which EETC is only a key part \cite{Li2014,Yang2015,Su2019,Particle1,Genetic1}. With the advancement of machine learning technology and availability of field data, some machine learning algorithms such as deep deterministic policy gradient \cite{Ning2021} and deep-neural-network based method \cite{Yin2020} were applied. Although the learning algorithms effectively reduce the dependency on the model construction, the global optimality cannot be guaranteed due to the existing local minimum \cite{Yin2020}. In addition, the computational time is still within a scale of seconds which could be a challenge for fast online operations \cite{Ning2021}. 

In recent years, an increasing number of studies apply mathematical programming to the EETC problem, which not only ensures the global optimality of the model but also gains high computational efficiency. Lu \textit{et al.} formulated the partial train speed trajectory optimization problem based on mixed-integer linear programming (MILP), which demonstrated the application potential of MILP for online computing \cite{Lu2016}. This method was later applied to the adaptive train speed trajectory optimization problem which can be extended to the EETC problem in real-world cases \cite{Tan2018}. Wang \textit{et al.} compared the performance of the pseudo-spectral method and MILP when dealing with the EETC problem, and concluded that the pseudo-spectral method has better optimality while MILP has higher computational efficiency \cite{Wang2013}. Researchers from Delft University of Technology conducted a series of studies based on the pseudo-spectral method to address the EETC problem and other relevant problems \cite{Wang2016,Wang2017a,Wang2020,Goverde2020}. Ye \textit{et al.} pointed out that piecewise linear (PWL) and piecewise quadratic functions can be used to approximate nonlinear engineering constraints and thus nonlinear programming can be applied to solve the EETC problem \cite{Ye2017}.

With many different methods studied on the field of EETC, convex optimization has rarely been applied. A problem can be solved efficiently if it can be identified or constructed as a convex optimization problem \cite{Convex_boyd}. Due to its fast-solvable characteristics, convex optimization has been applied in aerospace engineering for fast trajectory planning for Rendezvous and proximity operations \cite{Liu2013} and aerocapture \cite{Han2019}. Recently, researchers from the University of Birmingham proposed convex optimization to minimize the fuel consumption and concurrently optimizing the energy management of the hydrogen hybrid trains and applied the barrier method to solve the problem rapidly \cite{Jibrin2021,jibrin2021convex,jibrin2021mathematical}. In the series papers, EETC is not the only focus, though it serves as the fundamental result along with others such as the energy management strategy of the fuel cell hybrid train traction systems. The results demonstrated seem very promising and interesting to the field and the relaxation technology offers inspiration for us to bring the relaxation into the field of EETC. However, it remains an open question regarding the exactness of the optimization model. A much more thorough study is needed for exactness verification theoretically and experimentally.

In this paper, the EETC problem is formulated as a convex optimization model. Similar relaxation techniques proposed in \cite{Jibrin2021,jibrin2021convex,jibrin2021mathematical} have been applied to the nonconvex constraints regarding the calculation of time based on reciprocal of speed and kinetic energy based on the square of speed in the MILP-based model originally proposed in \cite{Tan2018}. The nonconvex constraints were initially addressed using the PWL modeling technique and now in this paper we are proposing to use convexificatin and relaxations. Accordingly, the MILP model can be converted into a convex optimization model and solved more efficiently.

This paper has the following contributions. 
\begin{itemize}
\item The exactness of the convex optimization model was verified using numerical experiments which indicate that the optimal solution of the proposed convex optimization model still appears in the feasible region of the original problem.
\item By comparing with the results of MILP and Radau pseudospectral method (RPM), the effectiveness of convex optimization was demonstrated, which also has a significant advantage in terms of computational efficiency. 
\item The computational advantage of the convex optimization model compared with MILP for large-scale problems was demonstrated by a series of numerical experiments with different model sizes.
\end{itemize} 

\section{Mathematical model}
\label{sec:model}
In this section, the EETC problem considering variable gradient and speed limit, punctuality, energy cost, and nonlinear motor characteristics was formulated initially. Two relaxations were applied to the nonconvex constraints subsequently, and the original model was converted to a convex model which can be solved efficiently.

\subsection{An energy-efficient train control model based on convex optimization}
\label{sec:EETC model}

\begin{figure}
    \centering
    \includegraphics{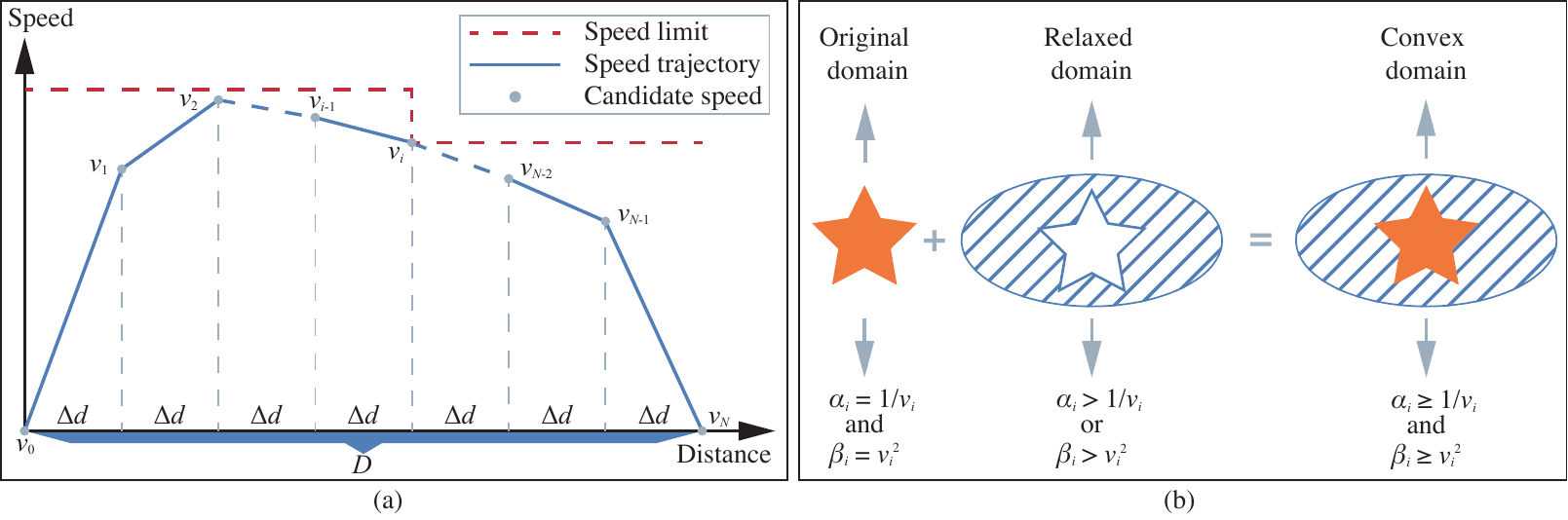}
    \caption{Schematic diagram of convex optimization for the energy-efficient train control problem. (a) Illustration of model construction of the mathematical programming. (b) The feasible domain change for the original problem after relaxations are applied. }
    \label{fig:modeling}
\end{figure}

In this paper, the running process of the train is constructed discretely, as shown in Figure~\ref{fig:modeling}(a). The whole journey is divided equally into $N$ segments with $\Delta d$, which can be expressed as
\begin{equation}
D = \sum_{i=1}^{N} \Delta d \label{con:delta_d}
\end{equation}
where $D$ is the total distance between two adjacent stations and $i$ represents the index of the discrete segments.

The whole journey includes $N+1$ candidate speed points $v_i$ corresponding to $N$ segments. Each speed point is supposed to be constrained by the speed limit $V_{i,max}$, which can be expressed as
\begin{equation}
0 < v_i \leq V_{i,max}.
\label{con:Speedlimit}
\end{equation}

Due to the punctuality requirement, the sum of elapsed time in each segment equals the total journey time $T$
\begin{equation}
\sum_{i=1}^{N} \frac{\Delta d}{v_i} = T.  \label{con:Time}
\end{equation}

During driving, the drag force always exists, which can be calculated by the Davis equation in \eqref{con:Davis}.
\begin{equation}
f_i = A + B v_{i} + C v_{i}^2 \label{con:Davis}
\end{equation}
where $A$, $B$ and $C$ are Davis coefficients.

When the train is in traction mode, mechanical energy is converted into kinetic energy, heat, and potential energy according to the conservation of energy. On the contrary, when the train conducts braking, the kinetic energy is converted into electric energy, heat, and potential energy. This process can be expressed as
\begin{equation}
F_i \Delta d - \frac{1}{2} M (v_{i}^2 - v_{i-1}^2) - f_i \Delta d - Mg \Delta H_i = 0 \label{con:Energy}
\end{equation}
where $F_i$ represents tractive (+) or braking (-) effort provided to the train, $M$ is the total mass of the train, $g$ is the gravitational constant, and $\Delta H_i$ is the altitude difference between the current position and previous position.

The electric energy consumed during traction and the regenerative energy recovered during braking can be calculated as
\begin{align}
	E_{i} &\geq F_i \Delta d / \eta_{t}\\
	E_{i} &\geq F_i \Delta d \eta_{b}
\end{align}
where $\eta_{t}$ is the traction efficiency of converting electrical energy into mechanical energy when $E_i > 0$ and $\eta_{b}$ is the braking efficiency of converting mechanical energy into electric energy when $E_i < 0$.

For the motor characteristics, inequality constraints \eqref{con:F} and \eqref{con:Power} are imposed to ensure that the train effort does not exceed the maximum tractive and braking effort.
\begin{equation}
-F_{b,max} \leq F_i \leq F_{t,max} \label{con:F}
\end{equation}
\begin{equation}
-P_{b,max}/v_i \leq F_i \leq P_{t,max}/v_i \label{con:Power}
\end{equation}

The net electrical energy of the entire journey is taken as the objective function, as shown in \eqref{con:Obj}.
\begin{equation}
\min \sum_{i=1}^{N} E_i \label{con:Obj}
\end{equation}

It can be found that the model proposed above is nonconvex due to the nonaffine equation \eqref{con:Time} and nonconvex inequality \eqref{con:Power} with $1/v_i$ and nonaffine equations \eqref{con:Davis} and \eqref{con:Energy} with $v_{i}^2$.

\subsection{Relaxations and convexification}
\label{sec:relaxations}

To convert the original nonconvex problem into a convex problem that can be solved efficiently by effective numerical methods like the interior-point method and barrier method, the following conditions need to be satisfied.
\begin{itemize}
\item The objective function must be convex.
\item The inequality constraint functions must be convex.
\item The equality constraint functions must be affine.
\end{itemize} 
To address the nonconvexity caused by the above-mentioned nonaffine equality and nonconvex inequality constraints, two relaxations were introduced to realize the convexification of the original problem, in which each nonconvex constraint is replaced with a relaxed but convex constraint. 

In the relaxation of \eqref{con:Time} and \eqref{con:Power}, $1/v_i$ is replaced by a auxiliary variable $\alpha_i$, as shown in \eqref{con:relaxT} and \eqref{con:relaxP}:
\begin{align}
\sum_{i=1}^{N} \Delta d \alpha_i&=T\label{con:relaxT}\\
-P_{b,max} \alpha_i \leq F_i &\leq P_{t,max} \alpha_i.
\label{con:relaxP}
\end{align}

Meanwhile, the inequality
\begin{equation}
\alpha_i \geq 1/v_i \label{con:Alpha}
\end{equation}
should be constructed to maintain the convexity of the model.

Likewise, $v_i^2$ is substituted by another auxiliary variable $\beta_i$ in the relaxation of \eqref{con:Davis} and \eqref{con:Energy} which can be rewritten as
\begin{equation}
f_i = A + B v_{i} + C \beta_i \label{con:RelaxDavis}
\end{equation}
\begin{equation}
F_i \Delta d - \frac{1}{2} M (\beta_{i} - \beta_{i-1}) - f_i \Delta d - Mg \Delta H_i = 0.
\label{con:RelaxEnergy}
\end{equation}

Inequality \eqref{con:Speedlimit} can be modified as
\begin{align}
\beta_i &\geq v_i^2  \label{con:Beta}\\
\beta_i &\leq V_{i,max}^2
\end{align}
to guarantee the convexity of the model while satisfying the speed limit.

In the paper \cite{Jibrin2021}, similar study has been conducted to optimize the traction energy without considering the characteristics of fuel cell and batteries in the proposed sequential method. The objective function thus takes the sum of traction energy plus $\gamma_i^2$ in each segment, where $\gamma_i$ is a variable similar to $\alpha_i$ in our proposed model which is different from our objective function but may not cause big variations given the small value of $\gamma_i$. In addition, our model takes both traction and regenerative braking energy into consideration which provides an extensions to applications involving regenerative braking.

By introducing auxiliary variables $\alpha_i$ and $\beta_i$, the EETC problem is transformed into a convex problem whose feasible domain is greater than that of the original problem, as shown in Figure~\ref{fig:modeling}(b). The orange star region indicates the original domain of the nonconvex EETC problem where relaxed constraints \eqref{con:Alpha} and \eqref{con:Beta} are held on the equal sign. The blue shaded region is generated by the relaxations and corresponds to the strict inequality relationship of relaxed constraints \eqref{con:Alpha} and \eqref{con:Beta}. The orange star region and the blue shaded region constitute the feasible region of the convex problem. It can be found that only the solutions from the orange star region in the convex domain can satisfy the physical constraints while solutions from the blue shaded region are not valid for the original EETC problem. Numerical experiments in Section~\ref{sec:simulation} conduct the verification of the model effectiveness of convex optimization, namely, verify that the optimal solution of the model appears in the orange region.

\section{Numerical experiments}
\label{sec:simulation}

In this section, numerical experiments were conducted to verify the effectiveness and fast solution performance of convex optimization for the EETC problem. Firstly, the effectiveness was illustrated indirectly by analyzing the optimal speed trajectories of the proposed model, MILP, and RPM. The modeling details of MILP and RPM can be found in references \cite{Tan2018,Goverde2020}. Secondly, the fast computation characteristic of convex optimization was demonstrated by comparing the CPU time of convex optimization with that of MILP and RPM with the same model size. Moreover, a comparative study with different model sizes was conducted to demonstrate the fast computational capability  of convex optimization in large-scale problems. Thirdly, the curves of the auxiliary variables and the corresponding original variables in a single train optimization were highlighted to verify the exactness of the proposed model.

In the optimization, route data between Expo Center Station and Wenquan East Station of Qingdao metro line 11 in China was applied to validate the proposed model. The rolling stock parameters are listed in Table~\ref{tab:Parameters of rolling stock} and the modeling parameters are shown in Table~\ref{tab:Parameters of three methods}. The PWL nodes used in the piecewise linearization of the MILP model are illustrated in Appendix~\ref{sec:PWL}. Parameter settings of GPOPS used to solve the RPM model are shown in Appendix~\ref{sec:RPM parameter}. Numerical experiments were carried out based on the computer configuration with 2.6 GHz \textit{Intel i7} CPU and 16 GB RAM.

\begin{table}
    \caption{Modelling parameters of rolling stock.}
    \label{tab:Parameters of rolling stock}
    \centering
    \begin{tabular}{p{8cm} p{2.5cm}<{\centering}}
    \toprule
    Parameter & Value \\
    \midrule
    $M (\rm t)$ & 144   \\
    $F_{t,max} (\rm kN)$ & 230.81    \\
    $F_{b,max} (\rm kN)$ & 230.81    \\
    $P_{t,max} (\rm kW)$ & 2520  \\
    $P_{b,max} (\rm kW)$ & 2520  \\
    $\eta_t$ & 0.9  \\
    $\eta_b$ & 0.6  \\
    $A ( \rm kN)$           & 3.0016    \\
    $B (\rm kN/(km/h))$     & 2.016e-2  \\
    $C (\rm kN/(km^2/h^2))$ & 6.9692e-4 \\
    \bottomrule
    \end{tabular}
\end{table}

\begin{table}
  \caption{Modeling parameters of the three methods.}
  \label{tab:Parameters of three methods}
  \centering
  \begin{tabular}{ccccc}
    \toprule
    Method & Solver & $D (\rm m)$ & $T (\rm s)$ & $N$\\
    \midrule
    Convex optimization  & GUROBI 9.5.0 & \multirow{3}{*}{2707}  & \multirow{3}{*}{140} & \multirow{3}{*}{238}\\
    MILP                 & GUROBI 9.5.0 &   &   &  \\
    RPM                  & GPOPS 4.1    &   &   &   \\
    \bottomrule
  \end{tabular}
\end{table}

\subsection{Effectiveness of convex optimization}

Note that regenerative braking was considered in the three methods. It should be declared that the whole journey in convex optimization and MILP was divided equally while the RPM model was divided into multiple segments with unequal distances. In spite of that, the total candidate speed points of the three methods were guaranteed to be consistent for the sake of a fair comparison.

\begin{figure}
    \centering
    \includegraphics{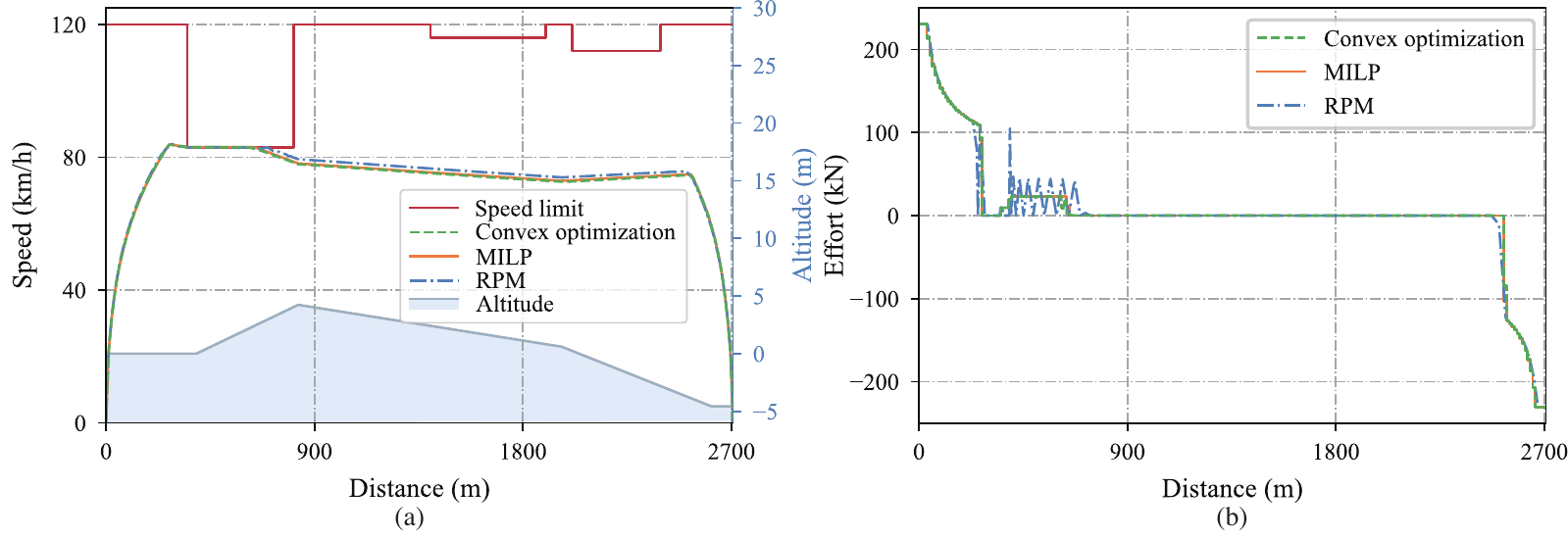}
    \caption{Optimal speed trajectory and tractive/braking efforts for real case based on convex optimization, mix-integer linear programming (MILP) and the Radau pseudospectral method (RPM).  (a) It illustrates the optimal speed trajectory with varying speed limit and gradient. (b) The tractive/braking effort shows that the train optimal operation consists of maximum traction, partial traction, coasting and maximum braking. }
    \label{fig:speed profile}
\end{figure}

The optimal speed trajectories and the tractive/braking efforts with the minimum energy consumption obtained by the three methods are shown in Figure~\ref{fig:speed profile}. As can be seen from Figure~\ref{fig:speed profile} (a), three speed trajectories indicate roughly the same trend in which the results of convex optimization and MILP are more smooth. RPM generates more fluctuations when conducting cruising which can be illustrated more explicitly in Figure~\ref{fig:speed profile} (b). Under the influence of gradient, the train speed changes during coasting. During the departure and arrival stages, the maximum traction and braking are adopted to accelerate and decelerate. 

The results including the objective function value, energy consumption, train target running time, train actual running time, and the CPU time of the optimization are shown in Table~\ref{tab:Results}. For the fact that the optimization results of the three methods include calculation errors introduced by the corresponding modeling characteristics, an independent simulation program was designed to recalculate the energy consumption and train running time with respect to each optimized speed trajectory. The pseudo code of recalculation program is shown in Appendix~\ref{sec:recalculation}. Thus, the energy consumption and actual train running time calculated independently were used to conduct the comparative analysis rather than the objective function values and the preset target journey time.

\begin{table}
 \centering
 \caption{Optimization results of the three EETC models}\label{tab:Results}
 \begin{threeparttable}
  \begin{tabular}{llll}
    \toprule
    Result & Convex optimization & MILP & RPM\\
    \midrule
    Objective function value (kWh)    & 9.335            & 9.404            & 9.458            \\
    Energy consumption\tnote{a} (kWh) & 9.330 (-0.005)   & 9.414 (+0.010)   & 9.869 (+0.411)   \\
    Target running time (s)           & 140              & 140              & 140              \\
    Actual running time\tnote{b} (s)  & 141.652 (+1.652) & 141.326 (+1.326) & 139.829 (-0.171) \\
    CPU time (s)                      & 0.160            & 0.900            & 5.290            \\
    \bottomrule
  \end{tabular}
  \begin{tablenotes}
  \footnotesize
    \item[a] The value in parenthesis is the difference between the simulated energy consumption and the objective function value.
    \item[b] The value in parenthesis is the difference between the actual train running time and the target journey time.
    \end{tablenotes}
    \end{threeparttable}
\end{table}

There is no significant difference in the total energy consumption obtained by the three methods. Due to the modeling error, the train running time of the RPM model is slightly shorter than the target time, while the total running time of the other two models is slightly longer than the target time, which leads to higher total energy consumption in the RPM model.

\subsection{Calculation performance of convex optimization}

In terms of CPU time, the computational efficiency of convex optimization is significantly higher than that of the other two mathematical models, as shown in Table~\ref{tab:Results}. It is 5.625 times over MILP and 33.063 times over RPM, which demonstrates the fast calculation capability of convex optimization and its potential in dealing with the online cases.

\begin{figure}
    \centering
    \includegraphics{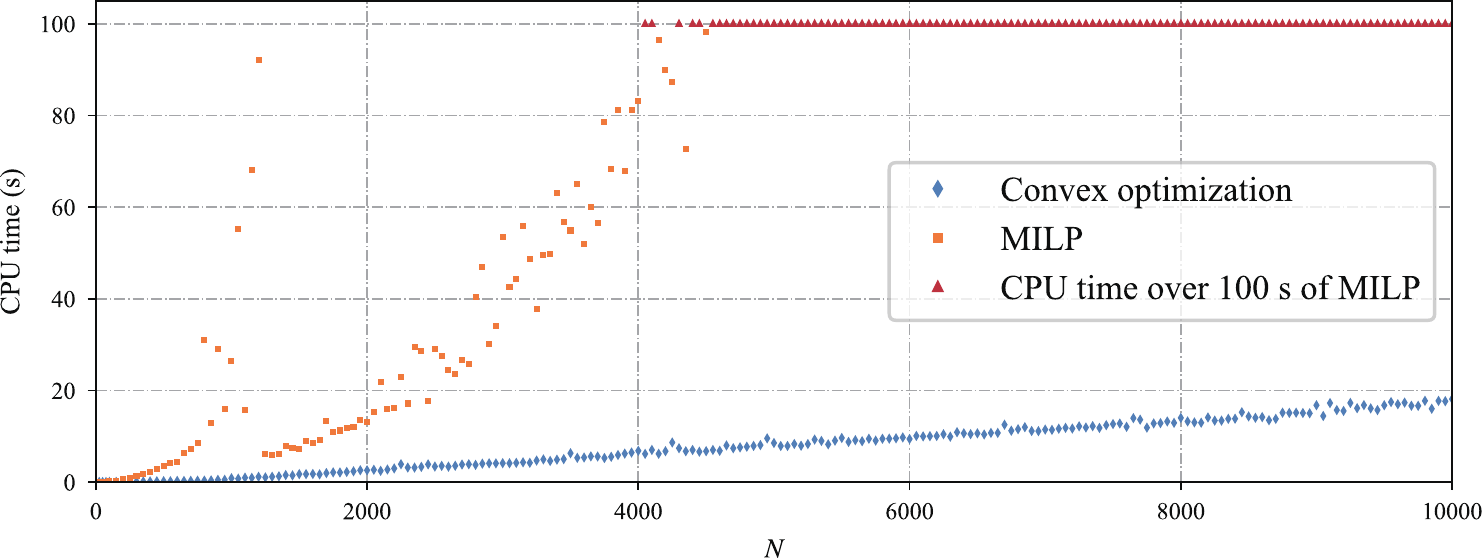}
    \caption{CPU time of convex optimization and MILP along with the segment number $N$. The green diamonds indicate that convex optimization takes no more than 20 s to obtain a solution as the model size increases. The orange squares illustrate the CPU time less than 100 s using MILP which has a faster growth trend than that of convex optimization. The red triangles show that in a series of optimizations, the CPU time of 56.93 \% of cases exceeded 100 s.}
    \label{fig:CPUtime}
\end{figure}

To figure out the calculation performance of convex optimization dealing with large-scale problems, multiple numerical experiments with variable segment number $N$ were conducted, as shown in Figure~\ref{fig:CPUtime}. In the series experiments, $N$ was set to be range from 25 to 10000 and the models being analyzed are the convex optimization model and MILP model. Due to the limitation caused by computational burden, we terminated and marked the experiments with a CPU time greater than 100 s. 

It can be observed from the result of convex optimization without any CPU time reaching 100 s, the maximum CPU time is 18.128 s with a segment number of 10000. The CPU time of 56.93 \% MILP optimization experiments exceeded 100 s, as indicated by the red triangles in Figure~\ref{fig:CPUtime}, most of which are large-scale optimizations ($N>4550$). In the MILP experiments with CPU time less than 100 s illustrated  orange squares, it can be found that the increase rate along with model size is faster than that of convex optimization. Overall, convex optimization demonstrates a faster computation performance when the scale of the optimization problem is greatly increased. 

\subsection{Validation of exactness}
\label{sec:validation of exactness}

\begin{figure}
    \centering
    \includegraphics{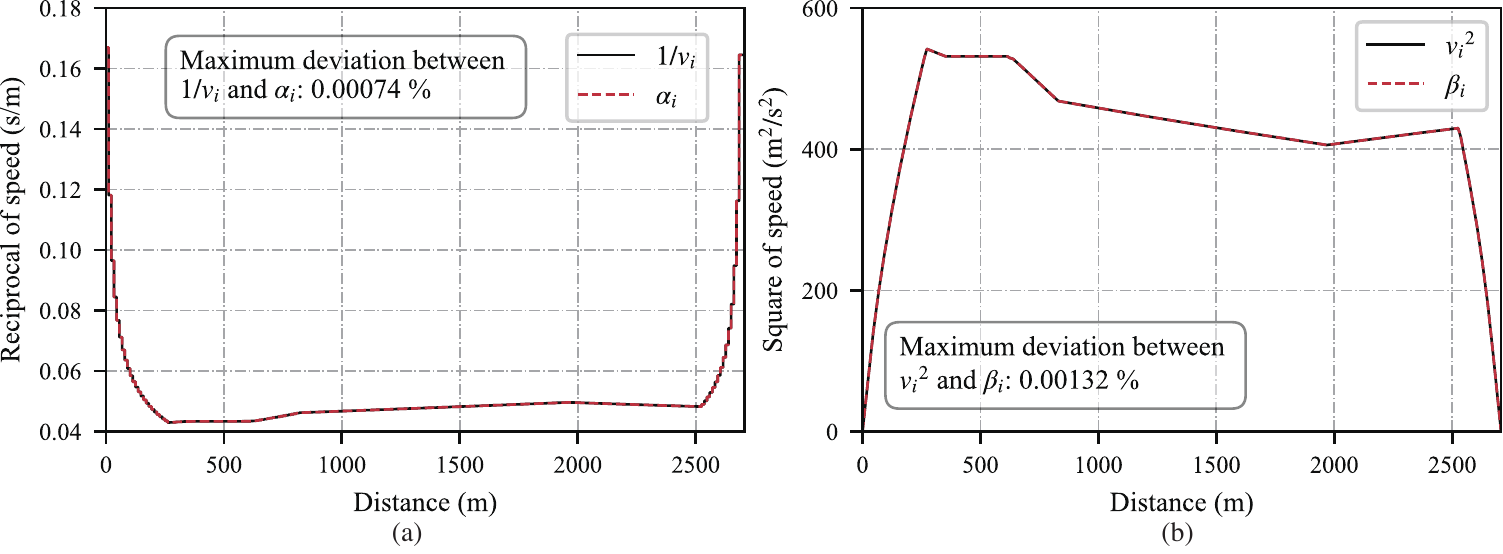}
    \caption{Comparison of values of the auxiliary variables and their original variables. (a) Relationship between the reciprocal of the candidate speed $1/v_i$ and the auxiliary variable $\alpha_i$. (b) Relationship between the square of the candidate speed $v_i^2$ and the auxiliary variable $\beta_i$. Both figures show that when the optimization model achieves the optimal solution, the value of the auxiliary variables are converging to the value of the original variables.}
    \label{fig:relaxation profile}
\end{figure}

Based on the above discussion in Section~\ref{sec:relaxations}, the optimal solution of the proposed model should be achieved with the relaxed inequality constraints set as equality, i.e. the auxiliary variables ($\alpha_i$ and $\beta_i$) are equal to the original variables ($1/v_i$ and $v_i^2$). To verify the exactness of the proposed model, Figure~\ref{fig:relaxation profile} depicts the values of $1/v_i$ and $v_i^2$ and those of the auxiliary variables $\alpha_i$ and $\beta_i$ respectively. The maximum deviations between the values of the initial variables and their auxiliary counterparts are 0.00074 \% ($1/v_i$ and $\alpha_i$) and 0.00132 \% ($v_i^2$ and $\beta_i$) respectively, which are negligible.

\section{Conclusion and future work}
In this paper, the energy-efficient train control (EETC) problem was formulated as a convex optimization problem and solved effectively and efficiently. Numerical experiments verified that the relaxed convex model still maintains exactness and the obtained solution is found in the feasible region of the original problem. The comparative study indicates that the solution of convex optimization is similar to that of mixed-integer linear programming (MILP) and the Radau pseudospectral method (RPM), but the computational efficiency is higher than the latter two mainstream methods often adopted by many other papers in the field.

The exactness verification in this paper is verified by numerical experiments, and the future work will be conducted to prove theoretically that when the convex optimization model attain the optimal solution at the equality of all the relaxation constraints.

\begin{appendices}
\newcounter{Afigure}
\setcounter{Afigure}{1}
\renewcommand\thefigure{\Alph{section}\arabic{Afigure}}
\newcounter{Atable}
\setcounter{Atable}{1}
\renewcommand\thetable{\Alph{section}\arabic{Atable}}

\section{Pisecewise linear nodes}\label{sec:PWL}

\begin{figure}
    \centering
    \includegraphics{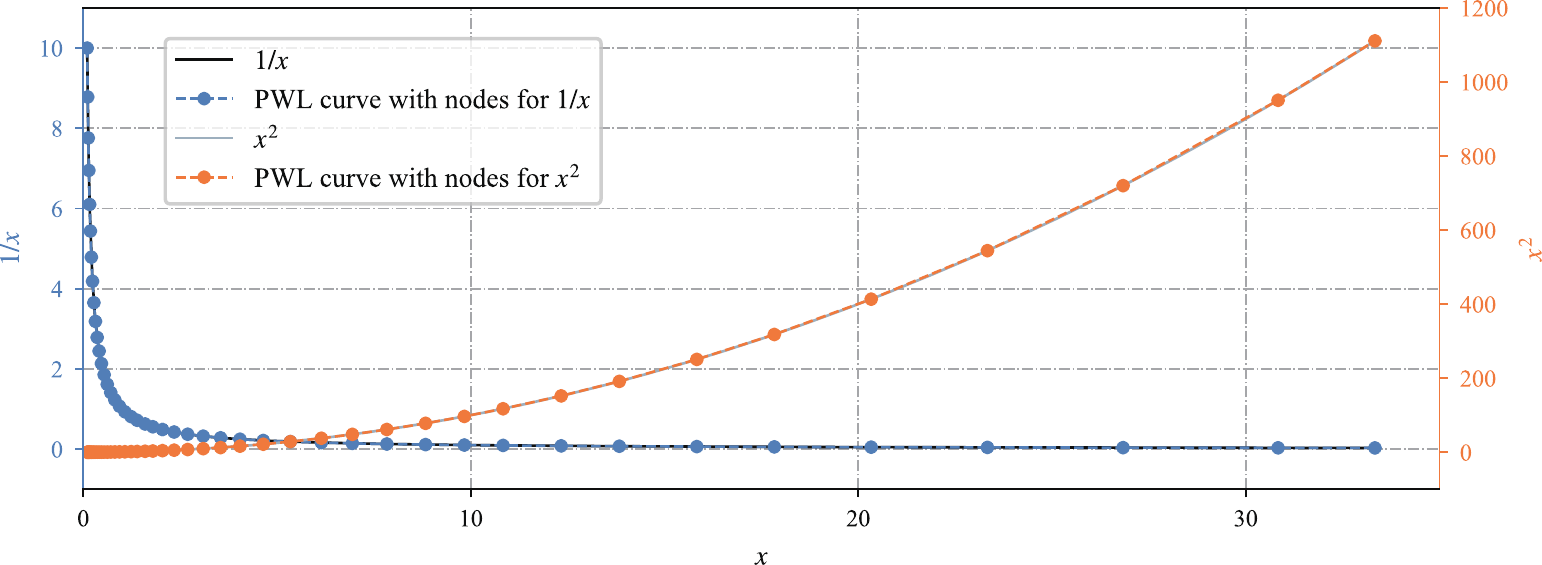}
    \caption{Piecewise linear nodes used in the mixed-integer linear programming for piecewise linearization.}
    \label{fig:appendix PWL}
\end{figure}

The selection of piecewise linear (PWL) nodes in mixed-integer linear programming (MILP) will affect the modeling accuracy and calculation efficiency of the model. We selected the number of nodes to ensure basic modeling accuracy. This appendix provides the PWL nodes selected in the paper, as illustrated in Figure~\ref{fig:appendix PWL} 

\section{Modeling parameters of Radau pseudospectral method}\label{sec:RPM parameter}

The solver settings adopted to solve the Radau pseudospectral method (RPM) model in this paper will affect the accuracy and computational efficiency of the optimization results, as shown in Table~\ref{tab:Parameters of GPOPS}.

\begin{table}

  \caption{Parameters settings of GPOPS Version 4.1}
  \label{tab:Parameters of GPOPS}
  \centering
  \begin{tabular}{ll}
    \toprule
    Parameter                  & Setting           \\
    \midrule
    setup.autoscale            & `on'              \\
    limits(p).nodesPerInterval & 25 (p=1,...,7)    \\
    limits(p).meshPoints       & [-1,1](p=1,...,7) \\
    setup.mesh.iteration       & 5                  \\
    setup.derivatives          & `finite-difference'          \\
    setup.tolerances           & 1e-1      \\
    \bottomrule
  \end{tabular}
\end{table}

\section{Verification simulation program for optimization results}\label{sec:recalculation}

Since the optimization results of the three methods analyzed in this paper include calculation errors introduced by the corresponding modeling characteristics, an independent verification simulation program was designed to recalculate the energy consumption and train running time with respect to a specific speed trajectory. The verification program pseudo code is shown in Procedure~\ref{alg1}. 

\floatname{algorithm}{Procedure}
\renewcommand{\algorithmicrequire}{\textbf{Input:}}
\renewcommand{\algorithmicensure}{\textbf{Output:}}
\begin{algorithm} 
	\caption{Verification program for the energy consumption $e$ and train running time $t$.} 
	\label{alg1} 
	\begin{algorithmic}[1]
	\REQUIRE A speed trajectory with $N+1$ candidate speed points ($d_i, v_i$).
	\STATE Initialization: $e \gets 0, t \gets 0$
    \FOR {$i = 1$; $i<N+1$; $i++$}
    \STATE $v_{ave} \gets (v_i+v_{i-1})/2$
    \STATE $\Delta d \gets d_i-d_{i-1}$
    \STATE $t \gets t + \Delta d/v_{ave}$
    \STATE Calculate energy consumed by drag:$e_1 \gets (A+Bv_{ave}+Cv_{ave}^2)\Delta d$
    \STATE Calculate kinetic energy change:$e_2 \gets 0.5M(v_i^2-v_{i-1}^2)$
    \STATE Calculate potential energy change:$e_3 \gets Mg\Delta H_i$
    \STATE Calculate mechanical energy of train:$e_4 \gets e_1 + e_2 + e_3$
    \IF{$e_4 \geq 0$}
    \STATE $e \gets e + e_4/\eta_t$
    \ELSE
    \STATE $e \gets e + e_4\eta_b$
    \ENDIF
    \ENDFOR
    \ENSURE $e, t$
	\end{algorithmic} 
\end{algorithm}
\end{appendices}

\bibliographystyle{unsrt}  
\bibliography{reference}
\end{document}